\documentclass{my_cmmse2012}

\usepackage[english]{babel}
\usepackage{hyperref}
\usepackage{amsmath}
\usepackage{amssymb}
\usepackage{amsthm}
\usepackage{graphicx}
\usepackage{newlfont}
\usepackage{setspace}
\usepackage{subfig}
\usepackage{color}
\usepackage{url}
\usepackage{fancybox} 
\usepackage[Lenny]{fncychap}
\usepackage{makeidx}
\usepackage{multirow}

\newtheorem{theorem}{Theorem}[section]

\begin{document}

\markboth{Modeling and Optimal Control}{H. S. Rodrigues, M. T. T. Monteiro, D. F. M. Torres}

\title{Modeling and Optimal Control\\
Applied to a Vector Borne Disease}


\author{Helena Sofia Rodrigues}{sofiarodrigues@esce.ipvc.pt}{1}

\author{M. Teresa T. Monteiro}{tm@dps.uminho.pt}{2}

\author{Delfim F. M. Torres}{delfim@ua.pt}{3}

\affiliation{1}{School of Business Studies}
  {Viana do Castelo Polytechnic Institute, Portugal}

\affiliation{2}{Department of Production and Systems}{University of Minho, Portugal}

\affiliation{3}{Department of Mathematics}{University of Aveiro, Portugal}


\begin{abstract}
A model with six mutually-exclusive compartments related to Dengue disease is presented. 
In this model there are three vector control tools: insecticides (larvicide and adulticide) 
and mechanical control. The problem is studied using an Optimal Control (OC) approach. 
The human data for the model is based on the Cape Verde Dengue outbreak. 
Some control measures are simulated and their consequences analyzed.


\msccodes 92B05, 93C95.

\end{abstract}


\section{Introduction}

Dengue is a vector borne disease transmitted to humans by the bite 
of an infected female \emph{Aedes} mosquito. Dengue transcends 
international borders and can be found in tropical and sub-tropical 
regions around the world, predominantly in urban and semi-urban areas. 
The risk may be aggravated further due to climate changes 
and to the globalization, as a consequence 
of the huge volume of international tourism and trade \cite{Semenza2009}.

There are four distinct, but closely related, viruses that cause
Dengue. Recovery from infection by one virus provides lifelong
immunity against that virus but confers only partial and transient
protection against subsequent infection by the other three
viruses \cite{Wearing2006}. Unfortunately, 
there is no specific effective treatment for Dengue.

Primary prevention of Dengue resides mainly in mosquito control.
There are two primary methods: larval control and adult mosquito
control, depending on the intended target \cite{Natal2002}. 
The application of adulticides can have a powerful impact 
on the abundance of adult mosquito vector. However, the efficacy 
is often constrained by the difficulty in achieving sufficiently 
high coverage of resting surfaces. This is the most common measure. 
Besides, the long term use of adulticide has several risks: 
the resistance of the mosquito to the product reducing 
its efficacy and the killing of other species that live 
in the same habitat. Larvicide treatment is done through 
long-lasting chemical in order to kill larvae and preferably 
have WHO clearance for use in drinking water \cite{Derouich2003}. 
Larvicide treatment is an effective way to control the vector
larvae, together with \emph{mechanical control}, which is related with
educational campaigns. The mechanical control must be done both 
by public health officials and by residents in affected areas. 
The participation of the entire population is essential 
to remove still water from domestic recipients, 
eliminating possible breeding sites \cite{Who2009}.

Mathematical modeling is an interesting tool for
understanding epidemiological diseases and for proposing 
effective strategies to fight them \cite{Lenhart2007}.


\section{Mathematical Model}

Taking into account the model presented in
\cite{Dumont2010,Dumont2008} and the considerations
of \cite{Sofia2009,Sofia2010c}, it is proposed
a new model more adapted
to the Dengue reality. The
notation used in the mathematical model includes three
epidemiological states for humans, indexed by $h$:
\begin{equation*}
\label{eq:star:1}
\begin{split}
& S_h \text{ susceptible (individuals who can contract the disease); }\\
& I_h \text{ infected (individuals capable of transmitting the disease to others); }\\
& R_h \text{ resistent (individuals who have acquired immunity). }
\end{split}
\end{equation*}

It is assumed that the total human population $(N_h)$
is constant and $N_h=S_h(t)+I_h(t)+R_h(t)$ at any time $t$.

There are three other state variables, related
to the female mosquitoes, indexed by $m$:
\begin{equation*}
\label{eq:star:2}
\begin{split}
& A_m \text{ aquatic phase (that includes the egg, larva and pupa stages);}\\
& S_m \text{ susceptible (mosquitoes that are able to contract the disease);}\\
& I_m \text{ infected (mosquitoes capable of transmitting the disease to humans).}
\end{split}
\end{equation*}

Due to short lifetime of mosquitoes (approximately 10 days), there is no resistant phase.
Humans and mosquitoes are assumed to be born susceptible.

To analyze the effect of campaigns in the combat
of the mosquito, three controls are considered:
\begin{equation*}
\label{eq:star:3}
\begin{split}
& c_{A} \text{ proportion of larvicide, and $0\leq c_A\leq1$;}\\
& c_{m} \text{ proportion of adulticide, and $0\leq c_m\leq1$;}\\
& \alpha \text{ proportion of mechanical control, and $0 < \alpha \leq1$.}
\end{split}
\end{equation*}

The aim of this work is to simulate different realities in order
to find the best policy to decrease the number of human infected.
A temporal mathematical model is introduced, with mutually-exclusive compartments,
to study the outbreak of 2009 in Cape Verde islands and
improving the model described in \cite{Sofia2009}.

The model considers the following parameters:
\begin{equation*}
\label{eq:star:4}
\begin{split}
& N_h \text{ total population;}\\
& B \text{ average daily biting (per day);}\\
& \beta_{mh} \text{ transmission probability from $I_m$ (per bite);}\\
& \beta_{hm} \text{ transmission probability from $I_h$ (per bite);}\\
& 1/\mu_{h} \text{ average lifespan of humans (in days);}\\
& 1/\eta_{h} \text{ mean viremic period (in days);}\\
& 1/\mu_{m} \text{ average lifespan of adult mosquitoes (in days);}\\
& \varphi \text{ number of eggs at each deposit per capita (per day);}\\
& 1/\mu_{A} \text{ natural mortality of larvae (per day);}\\
& \eta_{A} \text{ maturation rate from larvae to adult (per day);}\\
& m \text{ female mosquitoes per human;}\\
& k \text{ number of larvae per human.}
\end{split}
\end{equation*}

The Dengue epidemic is modeled by the following nonlinear
time-varying state equations:
\begin{equation}
\label{cap6_ode1}
\begin{cases}
\frac{dS_h}{dt}(t) = \mu_h N_h - \left(B\beta_{mh}\frac{I_m}{N_h}+\mu_h\right)S_h\\
\frac{dI_h}{dt}(t) = B\beta_{mh}\frac{I_m}{N_h}S_h -(\eta_h+\mu_h) I_h\\
\frac{dR_h}{dt}(t) = \eta_h I_h - \mu_h R_h\\
\frac{dA_m}{dt}(t) = \varphi \left(1-\frac{A_m}{\alpha k N_h}\right)(S_m+I_m)
-\left(\eta_A+\mu_A + c_A\right) A_m\\
\frac{dS_m}{dt}(t) = \eta_A A_m
-\left(B \beta_{hm}\frac{I_h}{N_h}+\mu_m + c_m\right) S_m\\
\frac{dI_m}{dt}(t) = B \beta_{hm}\frac{I_h}{N_h}S_m
-\left(\mu_m + c_m\right) I_m
\end{cases}
\end{equation}
with the initial conditions
\begin{equation}
\label{cap6_initial}
\begin{tabular}{llll}
$S_h(0)=S_{h0},$ &  $I_h(0)=I_{h0},$ &
$R_h(0)=R_{h0},$ \\
$A_m(0)=A_{m0},$ & $S_{m}(0)=S_{m0},$ & $I_m(0)=I_{m0}.$
\end{tabular}
\end{equation}


\section{Basic Reproduction Number}

An important measure of transmissibility of the disease 
is given by the basic reproduction number. It represents the
expected number of secondary cases produced in a completed
susceptible population, by a typical infected individual during
its entire period of infectiousness \cite{Hethcote2000}.
\medskip

\begin{theorem}
\label{thm:r0}
The basic reproduction number $\mathcal{R}_0$ associated to the
differential system \eqref{cap6_ode1} is given by
\begin{equation*}
\label{eq:R0}
\mathcal{R}_0 = \left(\frac{\alpha k B^2 \beta_{hm} \beta_{mh}
\mathcal{M}}{\varphi (\eta_h + \mu_h) (c_m + \mu_m)^2}\right)^{\frac{1}{2}}.
\end{equation*}
\end{theorem}

\medskip

\begin{proof}
The proof of this theorem is given in \cite{Sofia2013}.
\end{proof}

The model has two different populations (host and vector)
and the expected basic reproduction number should reflect the
infection transmitted from host to vector and vice-versa. 
If $\mathcal{R}_0<1$, then, on average, an infected individual
produces less than one new infected individual over the course of
its infectious period, and the disease cannot grow. Conversely, if
$\mathcal{R}_0>1$, then each individual infects more than one
person, and the disease can invade the population.

\begin{figure}
\centering
\subfloat[$\mathcal{R}_{0}$ as a function of $c_m$ and $c_A$]{\label{cap6_R0_cm_CA}
\includegraphics[width=0.5\textwidth]{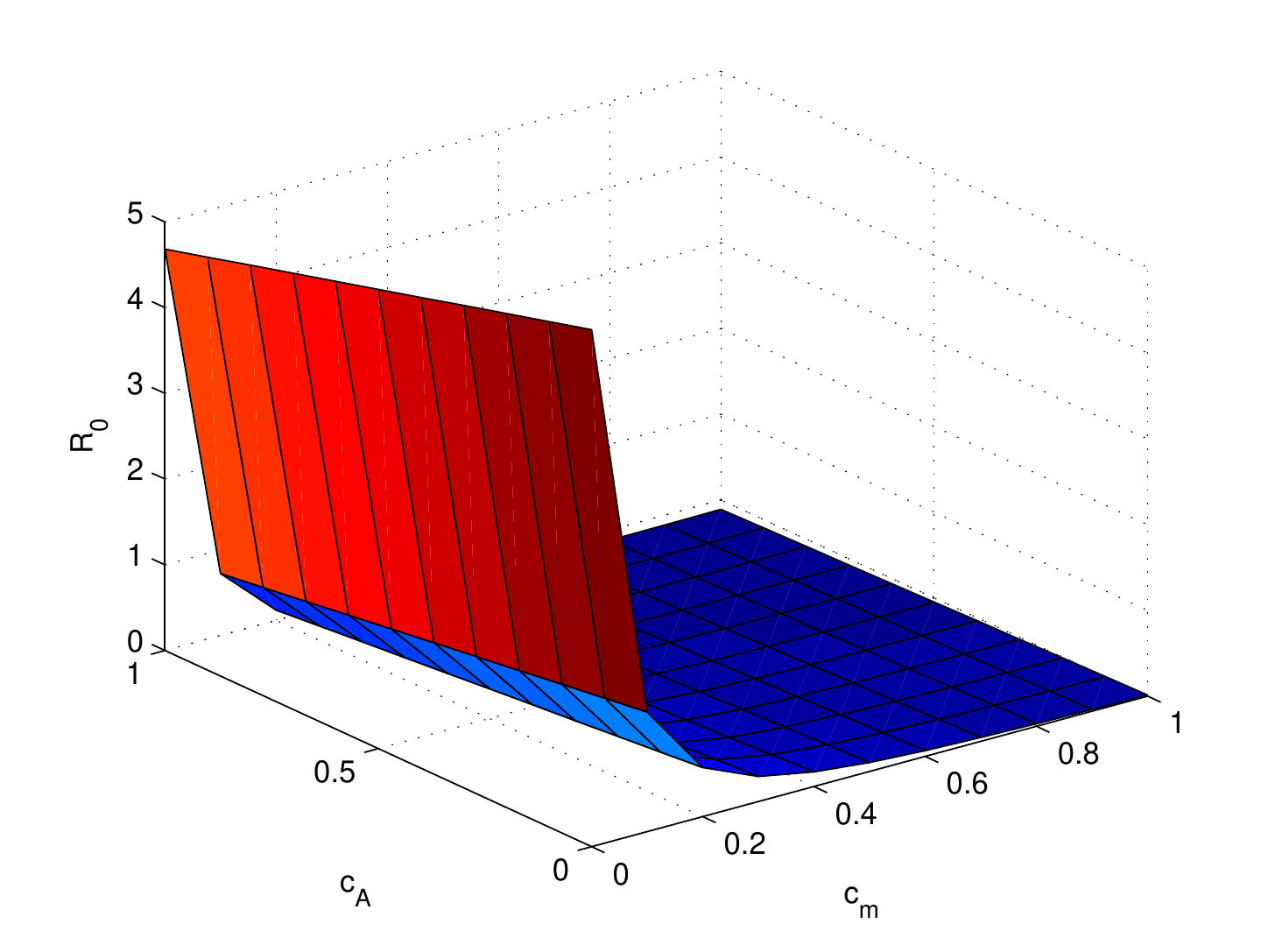}} ~ \\
\subfloat[$\mathcal{R}_{0}$ as a function of $c_m$ and $\alpha$]{
\label{cap6_R0_cm_alpha}
\includegraphics[width=0.5\textwidth]{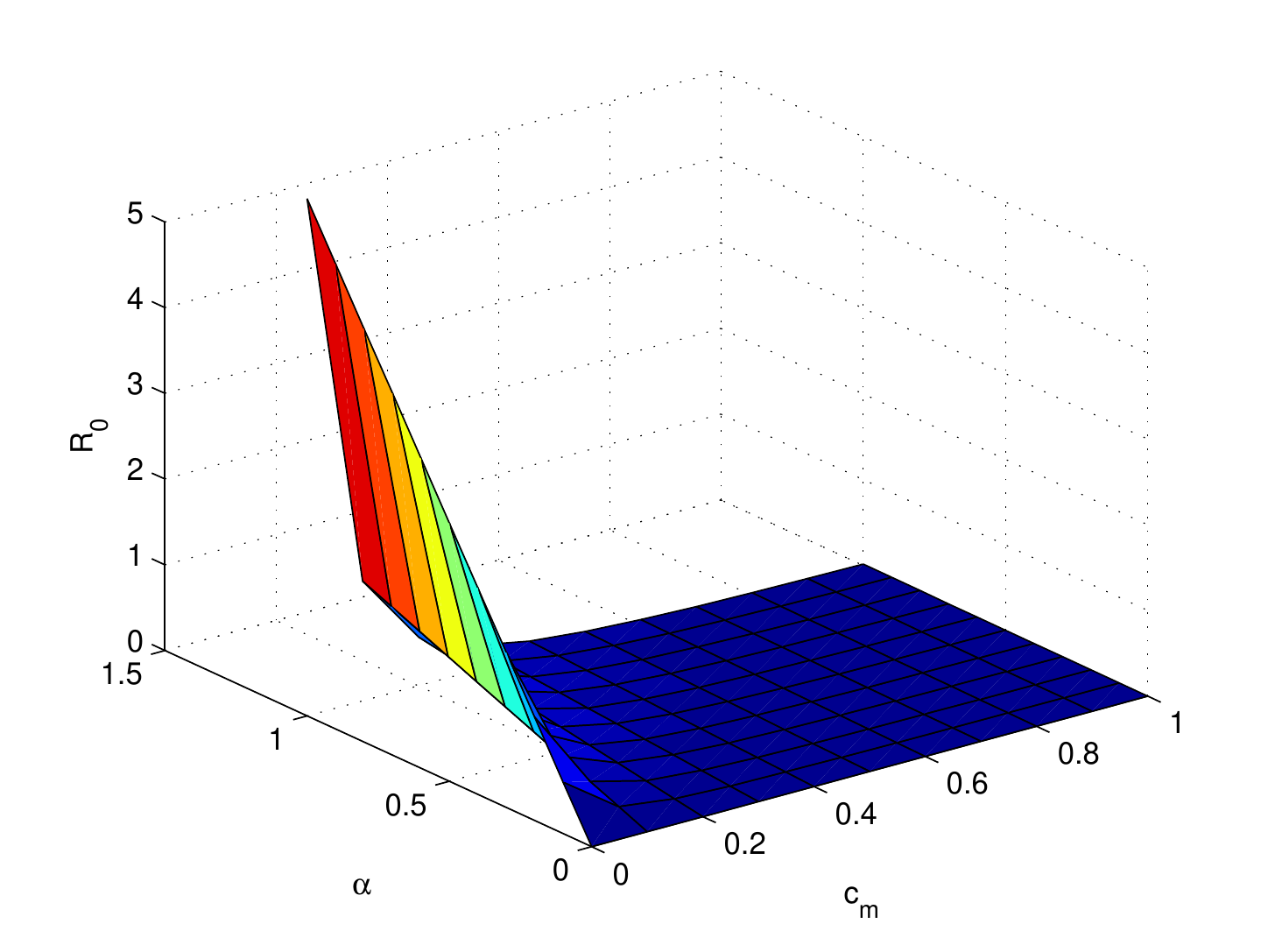}} ~ \\
\subfloat[$\mathcal{R}_{0}$ as a function of $c_A$ and $\alpha$]{
\label{cap6_R0_cA_alpha}
\includegraphics[width=0.5\textwidth]{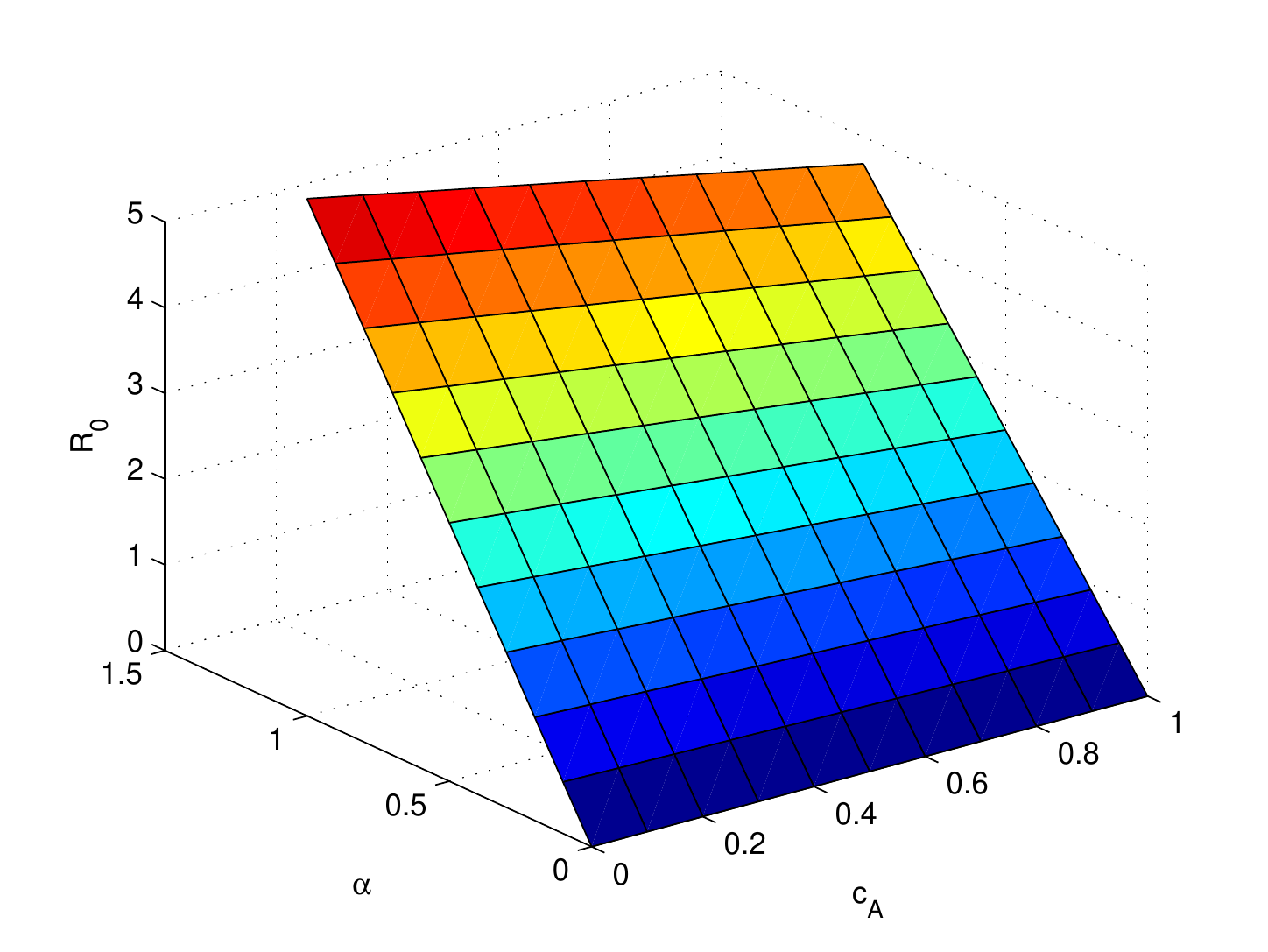}}
\caption{Influence of the controls on the basic reproduction number $\mathcal{R}_{0}$}
\label{cap6_R0_cm_CA_alpha}
\end{figure}

Assuming that parameters are fixed, the threshold $\mathcal{R}_0$ 
is influenceable by the control values. Figure \ref{cap6_R0_cm_CA_alpha} 
gives this relationship. It is possible to realize that the control $c_m$ 
is the one that most influences the basic reproduction number to stay below unit. 
Besides, the control in the aquatic phase alone is not enough to maintain 
$\mathcal{R}_0$ below unit: it requires an application close to 100\%.

\section{Optimal Control Approach}

Epidemiological models may give some basic guidelines for public health practitioners, 
comparing the effectiveness of different potential management strategies.

A cost functional was defined,

\begin{equation}
\label{chap6_functional}
\displaystyle J\left(u_1(.),u_2(.)\right)=\int_{0}^{t_f}\left[\gamma_D I_h(t)^2
+\gamma_S c_m(t)^2+\gamma_L c_A(t)^2+\gamma_E \left(1-\alpha\right)^2\right]dt,
\end{equation}
\noindent where $\gamma_D$, $\gamma_S$, $\gamma_L$ and $\gamma_E$ are weights related 
to the costs of the disease, adulticide, larvicide and mechanical control, 
respectively. In this way, an OC problem is defined:
\begin{equation*}
\begin{tabular}{ll}
minimize & (\ref{chap6_functional}) \\
subject to & (\ref{cap6_ode1}) and (\ref{cap6_initial}).
\end{tabular}
\end{equation*}


\section{Numerical Implementation}

The simulations were carried out using the following numerical
values: $N_h=480000$, $B=0.8$, $\beta_{mh}=0.375$,
$\beta_{hm}=0.375$, $\mu_{h}=1/(71\times365)$, $\eta_{h}=1/3$,
$\mu_{m}=1/10$,$\varphi=6$, $\mu_A=1/4$,
$\eta_A=0.08$, $m=3$, $k=3$. The initial conditions for the
problem were: $S_{h0}=N_h-10$, $I_{h0}=10$, $R_{h0}=0$,
$A_{m0}=k N_h$, $S_{m0}=m N_h$, $I_{m0}=0$.

In a first approach, the same weights were considered, 
which means $\gamma_D= \gamma_S=\gamma_L=\gamma_E=0.25$.

The OC problem was solved using two different softwares: DOTcvp and Muscod-II. 
The simulation behavior is similar, and we decide to show here only the results of the DOTcvp.

\begin{figure}[ptbh]
\begin{center}
\includegraphics[scale=0.7]{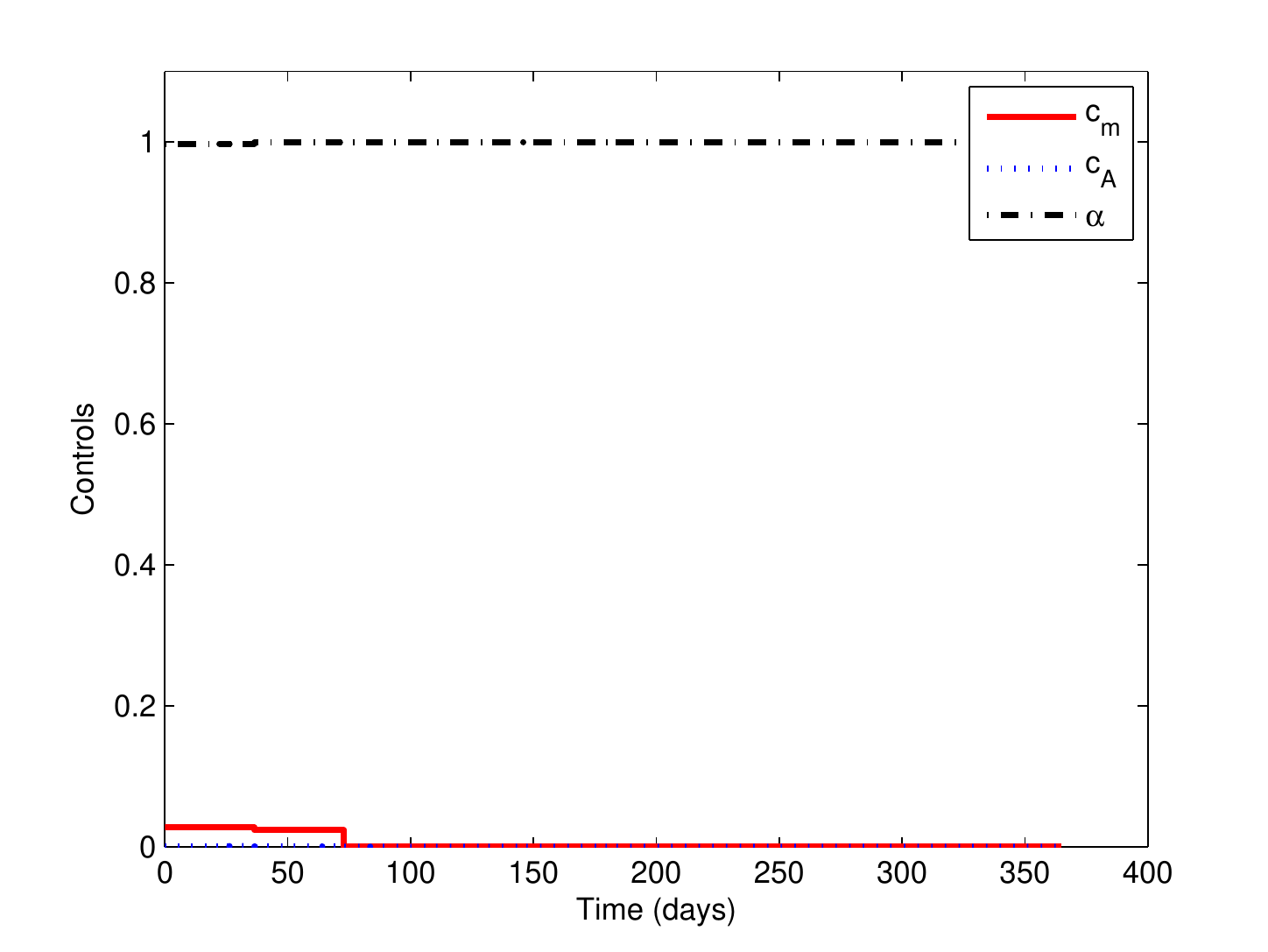}
\end{center}
\caption{Optimal control functions ($\gamma_D=\gamma_S=\gamma_L=\gamma_E=0.25$)}
\label{cap6_all_controls}
\end{figure}

The optimal functions for the controls are given in Figure~\ref{cap6_all_controls}.

As the results for the basic reproduction number, the adulticide was the control 
that more influences the decreasing of that ratio, and as consequence the decreasing 
of the number of infected persons and mosquitoes. Therefore, the adulticide 
was almost the one to be used. We believe  that the other controls do not assume 
an important role in the epidemic episode, because all the events happen 
in a short period of time, which means that adulticide has more impact. 
However the control of the mosquito in the aquatic phase can not be neglected. 
In situations of longer epidemic episodes or even in an endemic situation, 
the larval control represents an important tool.

\begin{figure}[ptbh]
\begin{center}
\includegraphics[scale=0.7]{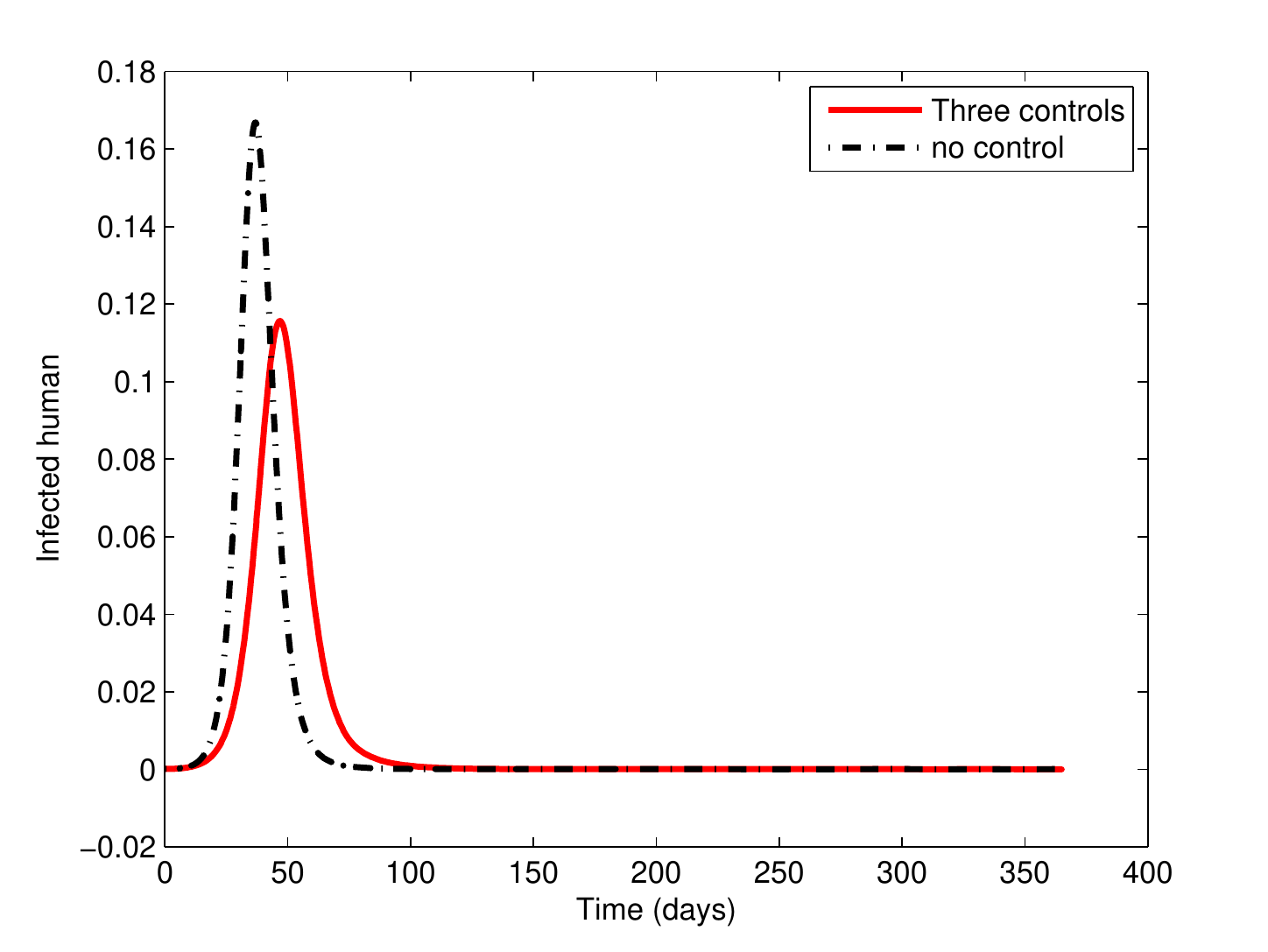}
\end{center}
\caption{Comparison of infected individuals under an optimal control situation and no controls.}
\label{cap6_all_controls_vs_nocontrol}
\end{figure}

Figure~\ref{cap6_all_controls_vs_nocontrol} presents the number of infected human. 
Comparing the optimal control case with a situation with no control, 
the number of infected people decreased considerably. Besides, in the situation 
where optimal control is used, the peak of infected people is minor, 
which facilitates the work in Health centers, because they can provide a better medical monitoring.

\section{Conclusions}
\label{sec:6:4}

A compartmental epidemiological model for Dengue disease was presented, 
composed by a set of differential equations. Simulations based on clean-up
campaigns to remove the vector breeding sites,
and also simulations on the application of insecticides
(larvicide and adulticide), were made. It was shown that even with
a low, although continuous, index of control over the time,
the results are surprisingly positive. The adulticide 
was the most effective control, from the fact that with a low percentage 
of insecticide, the basic reproduction number 
is kept below unit and the infected humans was smaller.

However, to bet only in adulticide is a risky decision. In some countries, 
such as Mexico and Brazil, the prolonged use of adulticides has been increasing 
the mosquito tolerance capacity to the product or even they become completely resistant. 
In countries where Dengue is a permanent threat, governments must act with differentiated tools. 
It will be interesting to analyze these controls in an endemic region 
and with several outbreaks. We claim that the results will be quite different.


\section*{Acknowledgements}

Work partially supported by the Portuguese Foundation for Science and Technology (FCT) through the Ph.D.
grant SFRH/BD/33384/2008 (Rodrigues) and the R\&D units Algoritmi (Monteiro) and CIDMA (Torres).

\small



\end{document}